\documentclass[review]{elsarticle}
\usepackage{graphicx}
\usepackage{epstopdf}
\usepackage{lineno,hyperref}
 \usepackage{geometry}
\usepackage{epsfig}
\usepackage{lineno,hyperref}
\usepackage{bookmark}
\usepackage{multirow}
 \usepackage{amsfonts}
\usepackage{multirow}
 \usepackage{algorithm}
 \usepackage[figuresright]{rotating}
 \usepackage{amscd}
\usepackage{mathrsfs}
\usepackage{booktabs,longtable}
\usepackage{indentfirst}
\usepackage{subfigure}
\usepackage{amstext}
\usepackage{amssymb}
\usepackage{fancyhdr}
\usepackage{amsmath}
 \usepackage{graphicx}
\usepackage{epstopdf}
\usepackage{lineno,hyperref}
 \usepackage{geometry}
\usepackage{epsfig}
\usepackage{lineno,hyperref}
\usepackage{bookmark}
\usepackage{multirow}
 \usepackage{amsfonts}
\usepackage{multirow}
 \usepackage{algorithm}
 \usepackage[figuresright]{rotating}
 \usepackage{amscd}
\usepackage{mathrsfs}
\usepackage{booktabs,longtable}
\usepackage{indentfirst}
\usepackage{subfigure}
\usepackage{amstext}
\usepackage{amssymb}
\usepackage{fancyhdr}
\usepackage{amsmath}

\modulolinenumbers[5]

\biboptions{sort&compress}
\newtheorem{thm}{Theorem}

\newdefinition{rmk}{Remark}
\newproof{pf}{Proof}
\newproof{pot}{Proof of Theorem \ref{thm2}}
\newtheorem{alg}{Algorithm}
\newdefinition{exa}{Example}

\journal{Journal of \LaTeX\ Templates}

\begin{document}

\begin{frontmatter}

\title{Convergence directions of the randomized Gauss--Seidel method and its extension \tnoteref{mytitlenote}}
\tnotetext[mytitlenote]{The work is supported by the National Natural Science Foundation of China (No. 11671060) and the Natural Science Foundation Project of CQ CSTC (No. cstc2019jcyj-msxmX0267)}

\author{Yanjun Zhang, \ Hanyu Li\corref{mycor}}
\cortext[mycor]{Corresponding author. E-mail addresses:  lihy.hy@gmail.com or hyli@cqu.edu.cn.}

\address{College of Mathematics and Statistics, Chongqing University, Chongqing 401331, P.R. China}

\begin{abstract}

The randomized Gauss--Seidel 
method and its extension have attracted much attention recently and  
their convergence rates have been considered extensively. However, the convergence rates are usually determined by upper bounds, which cannot fully reflect the  
actual convergence. In this paper, we make a detailed analysis of their convergence behaviors. The analysis shows that the larger the singular value of $A$ is, the faster the error decays in the corresponding singular vector space, and the convergence directions are mainly driven by the large singular values at the beginning, then gradually driven by the small singular values, and finally by the smallest nonzero singular value. These results explain the phenomenon found in the extensive numerical experiments appearing in the literature that these two methods seem to converge faster  
at the beginning. Numerical examples are provided to confirm the above findings.

\end{abstract}


\end{frontmatter}


\section{Introduction}
Linear least squares problem is a ubiquitous problem arising frequently in data analysis and scientific computing. Specifically, given a data matrix $A\in R^{m\times n}$ and a data vector $b\in R^{m}$, a linear least squares problem can be written as follows
\begin{equation}
\label{ls}
\min \limits _{ x \in R^{n}}\|b-Ax\|^2_{2}.
\end{equation}
In the literature, several direct methods have been proposed for solving its normal equations $A^TAx=A^Tb$ through either the QR factorization or the singular value decomposition (SVD) of $A^TA$ \cite{bjorck1996numerical, Higham2002}, which can be prohibitive when the matrix is large--scale. Hence, iterative methods are considered for solving large linear least squares problem, such as the famous Gauss--Seidel method \cite{Saad2003}.

In \cite{Leventhal2010}, Leventhal and Lewis proved that the randomized Gauss--Seidel (RGS) method, also known as the randomized coordinate descent method, converges to the solution at a linear rate in expectation. This method works on the columns of the matrix $A$ at random with probability proportional to their norms. Later, Ma, Needell and Ramdas \cite{Ma2015} provided a unified theory of the RGS method and the randomized Kaczmarz (RK) method \cite{Strohmer2009}, where the latter method works on the rows of $A$, and showed that the RGS method converges to the minimum Euclidean norm least squares solution $x_{\star}$ of (\ref{ls}) only when the matrix $A$ is of full column rank. To further develop the RGS method for more general matrix, inspired by the randomized extended Kaczmarz (REK) method \cite{Completion2013}, 
Ma et al. \cite{Ma2015} presented a variant of the RGS mehtod, i.e., randomized extended Gauss--Seidel (REGS) method, and proved that the REGS method converges to $x_{\star}$ regardless of whether the matrix $A$ has full column rank. After that, many variants of the RGS (or REGS) method were developed and studied extensively; see for example \cite{gower2015randomized, nutini2015coordinate, Hefny2017,tu2017breaking, xu2018hybrid,Dukui2019,razaviyayn2019linearly} and references therein.

To the best of our knowledge, when studying the convergence properties of the RGS and REGS methods, people mainly pay attention to their convergence rates and usually give corresponding upper bounds, and no 
work focuses on what determines their convergence rates, what drives their convergence directions, and what their ultimate directions is. 
As we know, the obtained upper bound of convergence can only be used as a reference for the convergence rate, and cannot truly reflect the empirical convergence of the method. So it is interesting to consider the above three problems.

In 2017, Jiao, Jin and Lu \cite{jiao2017preasymptotic} analyzed the preasymptotic convergence of the RK method. Recently, Steinerberger \cite{steinerberger2021randomized} made a more detailed analysis of the convergence property of the RK method for overdetermined full rank linear system. The author showed that the right singular vectors of the matrix $A$ describe the directions of distinguished dynamics and the RK method converges along small right singular vectors. 
After that, Zhang and Li \cite{zhang2021preconvergence} considered the convergence property of the REK method for all types of linear systems (consistent or inconsistent, overdetermined or underdetermined, full-rank or rank-deficient) and showed that the REK method converges to the minimum Euclidean norm least squares solution $x_{\star}$ with different decay rates in different right singular vectors spaces.

In this paper, we analyze the convergence properties of the RGS and REGS methods for linear least squares problem and show that the decay rates of the sequences $\{Ax_{k}\}_{k=1}^{\infty}$ and $\{ x_{k}\}_{k=1}^{\infty}$ (resp., the sequences $\{Az_{k}\}_{k=1}^{\infty}$ and $\{ z_{k}\}_{k=1}^{\infty}$ ) generated by the RGS method (resp., the REGS method) are depend on the size of singular values of $A$. Specifically, the larger the singular value of $A$ is, the faster the error decays in the corresponding singular vector space, and the convergence directions are mainly driven by the large singular values at the beginning, then gradually driven by the small singular values, and finally by the smallest nonzero singular value.

The rest of this paper is organized as follows. We first introduce some notations and preliminaries in Section \ref{sec2} and then present our main results about the RGS and REGS methods in Section \ref{sec3} and Section \ref{sec4}, respectively. Numerical experiments are given in Section \ref{sec5}.

\section{Notations and preliminaries }\label{sec2}
Throughout the paper, for a matrix $A$, $A^T$, $A^{(i)}$, $A_{(j)}$,
$\sigma_i(A)$, $\sigma_r(A)$, $\|A\|_F$, and $\mathcal{R}(A)$ denote its transpose, $i$th row (or $i$th entry in the case of a vector), $j$th column, $i$th singular value, smallest nonzero singular value, Frobenius norm, and column space, respectively. For any integer $m\geq1$, let $[m]:=\{1, 2, 3, ..., m\}$. If the matrix $G\in R^{n\times n}$ is positive definite, we define the energy norm of any vector $x\in R^{n}$ as $\| x\|_G:=\sqrt{x^TGx}$. In addition, we denote the identity matrix by $I$, its $j$th column by $e_{(j)}$ and the expectation of any random variable $\xi$ by $\mathbb{E} [\xi]$.

In the following, we use $x_{\star}=A^{\dag}b$ to denote the minimum Euclidean norm least squares solution of (\ref{ls}), where $A^{\dag}$ denotes the Moore--Penrose pseudoinverse of the matrix $A$. Because the SVD is the basic tool for the convergence analysis in next two sections, we denote the SVD \cite{golub2013matrix} of $A\in R^{m\times n}$ by
\begin{align}
A=U\Sigma V^{T}, \notag
\end{align}
where $U=[u_1, u_2, \ldots u_m]\in R^{m\times m}$ and $V=[v_1, v_2, \ldots v_n]\in R^{n\times n}$ are column orthonormal matrices and their column vectors known as the left and right singular vectors, respectively, and $\Sigma\in R^{m\times n}$ is diagonal with the diagonal elements ordered nonincreasingly, i.e., $\sigma_1(A)\geq \sigma_2(A)\geq \ldots \sigma_r(A)>0$ with $r\leq \min\{m, n\}$.

\section{Convergence directions of the RGS method}\label{sec3}
We first list the RGS method \cite{Leventhal2010, Ma2015} in Algorithm \ref{alg1} and restate its convergence bound in Theorem \ref{theorem0}.
\begin{alg}
\label{alg1}
 The RGS method
\begin{enumerate}[]
\item \mbox{INPUT:} ~$A$, $b$, $\ell$, $x_{0}\in R^{n }$
\item For $k=1, 2, \ldots, \ell-1$ do
\item ~~~~Select $j\in [n]$ with probability $\frac{\|A_{(j)}\|^2_2}{\|A\|^2_F}$
\item ~~~~Set $x_{k}=x_{k-1}-\frac{A_{(j)}^T ( Ax_{k-1}-b)}{ \| A_{(j)} \|_{2}^{2}}e_{(j)}$
\item End for
\end{enumerate}
\end{alg}

\begin{thm}
\label{theorem0}\cite{Leventhal2010, Ma2015}
Let $A\in R^{m\times n}$, $b\in R^{m}$, $x_{\star}=A^{\dag}b$ be the minimum Euclidean norm least squares solution, and $x_k$ be the $k$th approximation of the RGS method generated by Algorithm \ref{alg1} with initial guess $x_{0}\in R^{n }$. Then
\begin{align}
\mathbb{E}[\| Ax_{k}- Ax_{\star}\|_2^2]\leq(1-\frac{\sigma_r^2(A)}{\|A\|^2_F})^k \|Ax_{0}- Ax_{\star}\|_2^2.\label{th0}
\end{align}
\end{thm}

\begin{rmk}
\label{rmk1}
 Theorem \ref{theorem0} shows that $Ax_{k}$ converges linearly in expectation to $Ax_{\star}$ regardless of whether the matrix $A$ has full rank. Since $\| Ax_{k}- Ax_{\star}\|_2^2=\|  x_{k}- x_{\star}\|_{A^TA}^2$, it follows from (\ref{th0}) that
\begin{align}
\mathbb{E}[\| x_{k}-  x_{\star}\|_{A^TA}^2]\leq(1-\frac{\sigma_r^2(A)}{\|A\|^2_F})^k \| x_{0}-  x_{\star}\|_{A^TA}^2,\notag
\end{align}
which implies that $ x_{k}$ converges linearly in expectation to the minimum Euclidean norm least squares solution $ x_{\star}$ when the matrix $A$ is overdetermined and of full column rank, but can not converge to $ x_{\star}$ when $A$ is not full column rank. So, we 
assume that the matrix $A$ is of full column rank in this section.
\end{rmk}

Now, we give our three main results of the RGS method.
\begin{thm}
\label{theorem1}
Let $A\in R^{m\times n}$, $b\in R^{m}$, $x_{\star}=A^{\dag}b$ be the minimum Euclidean norm least squares solution, and $x_k$ be the $k$th approximation of the RGS method generated by Algorithm \ref{alg1} with initial guess $x_{0}\in R^{n }$. Then
\begin{align}
\mathbb{E}[\langle Ax_{k}- Ax_{\star}, u_{\ell} \rangle]= (1-\frac{\sigma_\ell^2(A)}{\|A\|^2_F})^k \langle Ax_{0}- Ax_{\star}, u_{\ell} \rangle.\label{th1}
\end{align}
\end{thm}

\begin{pf}
Let $\mathbb{E}_{k-1}[\cdot]$ be the conditional expectation conditioned on the first $k-1$ iterations of the RGS method. Then, from Algorithm \ref{alg1}, we have
\begin{align}
 & \mathbb{E}_{k-1}[\langle Ax_{k}- Ax_{\star}, u_{\ell} \rangle]\notag
\\
 &= \sum\limits_{j=1}^{n}\frac{ \|A_{ (j )} \|_{2}^{2}}{\|A\|_{F}^{2}} \langle Ax_{k-1}-\frac{A_{ (j )}^T(Ax_{k-1}-b)}{\|A_{ (j )} \|_{2}^{2}}A_{ (j )}- Ax_{\star}, u_{\ell} \rangle \notag
\\
 &=  \langle Ax_{k-1}- Ax_{\star}, u_{\ell} \rangle - \frac{ 1}{\|A\|_{F}^{2}} \sum\limits_{j=1}^{n} \langle A_{ (j )}^T(Ax_{k-1}-b) A_{ (j )} , u_{\ell}  \rangle  \notag
\\
 &=  \langle Ax_{k-1}- Ax_{\star}, u_{\ell} \rangle - \frac{ 1}{\|A\|_{F}^{2}} \sum\limits_{j=1}^{n} \langle A_{ (j )},Ax_{k-1}-b  \rangle \langle A_{ (j )}, u_{\ell}  \rangle  \notag
\\
 &=  \langle Ax_{k-1}- Ax_{\star}, u_{\ell} \rangle - \frac{ 1}{\|A\|_{F}^{2}}  \langle A^T(Ax_{k-1}-b)   , A^Tu_{\ell}  \rangle,  \notag
\end{align}
which together with the facts $A^T(b-Ax_{\star})=0$ and $A^Tu_{\ell}=\sigma_{\ell}(A) v_{\ell} $ yields
\begin{align}
 & \mathbb{E}_{k-1}[\langle Ax_{k}- Ax_{\star}, u_{\ell} \rangle]\notag
\\
 &=  \langle Ax_{k-1}- Ax_{\star}, u_{\ell} \rangle - \frac{ 1}{\|A\|_{F}^{2}}  \langle A^T(Ax_{k-1}-Ax_{\star})   , A^Tu_{\ell}  \rangle  \notag
\\
 &=  \langle Ax_{k-1}- Ax_{\star}, u_{\ell} \rangle - \frac{ 1}{\|A\|_{F}^{2}}  \langle A^T(\sum\limits_{i=1}^{m} \langle Ax_{k-1}-Ax_{\star}, u_i\rangle u_i)   , A^Tu_{\ell}  \rangle   \notag
\\
 &=  \langle Ax_{k-1}- Ax_{\star}, u_{\ell} \rangle - \frac{ 1}{\|A\|_{F}^{2}}  \langle  (\sum\limits_{i=1}^{m} \langle Ax_{k-1}-Ax_{\star}, u_i\rangle \sigma_{i}(A) v_{i})   , \sigma_{\ell}(A) v_{\ell}  \rangle   \notag
\\
 &=  \langle Ax_{k-1}- Ax_{\star}, u_{\ell} \rangle - \frac{\sigma_{\ell}^2(A)}{\|A\|_{F}^{2}}       \langle Ax_{k-1}-Ax_{\star}, u_{\ell}\rangle    \notag
\\
 &=  (1- \frac{\sigma_{\ell}^2(A)}{\|A\|_{F}^{2}}   )   \langle Ax_{k-1}-Ax_{\star}, u_{\ell}\rangle.    \notag
\end{align}
Thus, by taking the full expectation on both sides, we have
\begin{align}
 \mathbb{E}[\langle Ax_{k}- Ax_{\star}, u_{\ell} \rangle] =  (1- \frac{\sigma_{\ell}^2(A)}{\|A\|_{F}^{2}}   )  \mathbb{E}[ \langle Ax_{k-1}-Ax_{\star}, u_{\ell}\rangle  ] = \ldots= (1- \frac{\sigma_{\ell}^2(A)}{\|A\|_{F}^{2}}   )^k   \langle Ax_{0}-Ax_{\star}, u_{\ell}\rangle,   \notag
\end{align}
which is the estimate (\ref{th1}).

\end{pf}

\begin{rmk}
\label{rmk2}
 Theorem \ref{theorem1} shows that the decay rates of $\|Ax_k-Ax_{\star}\|_2$ are different in different left singular vectors spaces. Specifically, the decay rates are dependent on the singular values: the larger the singular value of $A$ is, the faster the error decays in the corresponding left singular vector space. This implies that the smallest singular value will lead to the slowest rate of convergence, which is the one in (\ref{th0}). So, the convergence bound presented in \cite{Leventhal2010, Ma2015} is optimal.
\end{rmk}

\begin{rmk}
\label{rmk3}
Let $r_k=b-Ax_k$ be the residual vector with respect to the $k$-th approximation $x_k$, and $r_{\star}=b-Ax_{\star}$ be the true residual vector with respect to the minimum Euclidean norm least squares solution $x_{\star}$. It follows from (\ref{th1}) and $Ax_k-Ax_{\star}=-( r_{k}- r_{\star})$ that
\begin{align}
\mathbb{E}[\langle r_{k}- r_{\star}, u_{\ell} \rangle]= (1-\frac{\sigma_\ell^2(A)}{\|A\|^2_F})^k \langle r_{0}- r_{\star}, u_{\ell} \rangle.\notag
\end{align}
Hence, Theorem \ref{theorem1} also implies that the decay rates of $ \| r_{k}- r_{\star}\|_2 $ of the RGS method depend on the singular values.
\end{rmk}

\begin{rmk}
\label{rmk4}
 Using the facts $\langle Ax_{k}- Ax_{\star}, u_{\ell} \rangle=\langle  x_{k}- x_{\star},A^Tu_{\ell} \rangle$ and $A^Tu_{\ell}=\sigma_{\ell}(A) v_{\ell} $, from (\ref{th1}), we have
 \begin{align}
\mathbb{E}[\langle  x_{k}-  x_{\star},  v_{\ell}  \rangle]= (1-\frac{\sigma_\ell^2(A)}{\|A\|^2_F})^k \langle  x_{0}-  x_{\star},   v_{\ell}  \rangle,\notag
\end{align}
which recovers the decay rates of the RK method in different right singular vectors spaces \cite{steinerberger2021randomized}. In this view, both RGS and RK methods are essentially equivalent.
\end{rmk}

\begin{thm}
\label{theorem2}
Let $A\in R^{m\times n}$, $b\in R^{m}$, $x_{\star}=A^{\dag}b$ be the minimum Euclidean norm least squares solution, and $x_k$ be the $k$th approximation of the RGS method generated by Algorithm \ref{alg1} with initial guess $x_{0}\in R^{n }$. Then
\begin{align}
\mathbb{E}[\|  Ax_{k}- Ax_{\star}\|_2^2]= \mathbb{E}[(1-\frac{1}{\|A\|^2_F}\|A^T\frac{Ax_{k-1}- Ax_{\star}}{\|Ax_{k-1}- Ax_{\star}\|_2}\|_2^2)\|  Ax_{k-1}- Ax_{\star}\|_2^2]. \notag
\end{align}
\end{thm}

\begin{pf}
Similar to the proof of \cite{Ma2015}, we can derive the desired result.
\end{pf}

\begin{rmk}
\label{rmk5}
Since $\|A^T\frac{Ax_{k-1}- Ax_{\star}}{\|Ax_{k-1}- Ax_{\star}\|_2}\|_2^2\geq\sigma_r^2(A)$, Theorem \ref{theorem2} implies that the RGS method actually converges faster if $Ax_{k-1}-Ax_{\star}$ is not close to left singular vectors corresponding to the small singular values of $A$ .
\end{rmk}

\begin{thm}
\label{theorem3}
Let $A\in R^{m\times n}$, $b\in R^{m}$, $x_{\star}=A^{\dag}b$ be the minimum Euclidean norm least squares solution, and $x_k$ be the $k$th approximation of the RGS method generated by Algorithm \ref{alg1} with initial guess $x_{0}\in R^{n }$. Then
\begin{align}
\mathbb{E}[\langle  \frac{Ax_{k}- Ax_{\star}}{\|Ax_{k}- Ax_{\star}\|_2}, \frac{Ax_{k+1}- Ax_{\star}}{\|Ax_{k+1}- Ax_{\star}\|_2}  \rangle^2]= 1 -\frac{1}{\|A\|^2_F} \mathbb{E}[ \|A^T\frac{Ax_{k }- Ax_{\star}}{\|Ax_{k }- Ax_{\star}\|_2}\|_2^2  ]. \label{th3}
\end{align}
\end{thm}

\begin{pf}
From Algorithm \ref{alg1}, we have
\begin{align}
 & \mathbb{E}_{k}[\langle  \frac{Ax_{k}- Ax_{\star}}{\|Ax_{k}- Ax_{\star}\|_2}, \frac{Ax_{k+1}- Ax_{\star}}{\|Ax_{k+1}- Ax_{\star}\|_2}  \rangle^2]\notag
\\
 &= \mathbb{E}_{k}[\langle  \frac{Ax_{k}- Ax_{\star}}{\|Ax_{k}- Ax_{\star}\|_2}, \frac{Ax_{k}-\frac{A_{ (j )}^T(Ax_{k}-b)}{\|A_{ (j )} \|_{2}^{2}}A_{ (j )}- Ax_{\star}}{\|Ax_{k}-\frac{A_{ (j )}^T(Ax_{k}-b)}{\|A_{ (j )} \|_{2}^{2}}A_{ (j )}- Ax_{\star}\|_2}  \rangle^2]\notag
\\
 &= \mathbb{E}_{k}[\frac{1}{\|Ax_{k}- Ax_{\star}\|_2^2\cdot \|Ax_{k}-\frac{A_{ (j )}^T(Ax_{k}-b)}{\|A_{ (j )} \|_{2}^{2}}A_{ (j )}- Ax_{\star}\|_2^2}\langle   Ax_{k}- Ax_{\star} ,  Ax_{k}-\frac{A_{ (j )}^T(Ax_{k}-b)}{\|A_{ (j )} \|_{2}^{2}}A_{ (j )}- Ax_{\star} \rangle^2]. \notag
\end{align}
Since $\langle   Ax_{k}- Ax_{\star} ,  Ax_{k}-\frac{A_{ (j )}^T(Ax_{k}-b)}{\|A_{ (j )} \|_{2}^{2}}A_{ (j )}- Ax_{\star} \rangle=\|Ax_{k}-\frac{A_{ (j )}^T(Ax_{k}-b)}{\|A_{ (j )} \|_{2}^{2}}A_{ (j )}- Ax_{\star}\|_2^2$, we have
\begin{align}
 & \mathbb{E}_{k}[\langle  \frac{Ax_{k}- Ax_{\star}}{\|Ax_{k}- Ax_{\star}\|_2}, \frac{Ax_{k+1}- Ax_{\star}}{\|Ax_{k+1}- Ax_{\star}\|_2}  \rangle^2]\notag
\\
 &= \mathbb{E}_{k}[\frac{ 1}{\|Ax_{k}- Ax_{\star}\|_2^2 }\|Ax_{k}-\frac{A_{ (j )}^T(Ax_{k}-b)}{\|A_{ (j )} \|_{2}^{2}}A_{ (j )}- Ax_{\star}\|_2^2 ] \notag
\\
 &= \mathbb{E}_{k}[\frac{ 1}{\|Ax_{k}- Ax_{\star}\|_2^2 }(\|Ax_{k} - Ax_{\star}\|_2^2-2 \langle Ax_{k}- Ax_{\star}, \frac{A_{ (j )}^T(Ax_{k}-b)}{\|A_{ (j )} \|_{2}^{2}}A_{ (j )} \rangle +\frac{(A_{ (j )}^T(Ax_{k}-b))^2}{\|A_{ (j )} \|_{2}^{2}})  ] \notag
\\
 &= \mathbb{E}_{k}[\frac{ 1}{\|Ax_{k}- Ax_{\star}\|_2^2 }(\|Ax_{k} - Ax_{\star}\|_2^2-\frac{(A_{ (j )}^T(Ax_{k}-Ax_{\star}))^2}{\|A_{ (j )} \|_{2}^{2}})  ] \notag
\\
 &= \mathbb{E}_{k}[1-\frac{(A_{ (j )}^T  \frac{Ax_{k}-Ax_{\star}}{\|Ax_{k}- Ax_{\star}\|_2  } )^2}{\|A_{ (j )} \|_{2}^{2}}  ] \notag
\\
 &= \sum\limits_{j=1}^{n}\frac{ \|A_{ (j )} \|_{2}^{2}}{\|A\|_{F}^{2}} (1-\frac{(A_{ (j )}^T  \frac{Ax_{k}-Ax_{\star}}{\|Ax_{k}- Ax_{\star}\|_2  } )^2}{\|A_{ (j )} \|_{2}^{2}})   \notag
\\
 &=  1-\frac{ 1}{\|A\|_{F}^{2}} \| A^T  \frac{Ax_{k}-Ax_{\star}}{\|Ax_{k}- Ax_{\star}\|_2  }  \|_2^2.  \notag
\end{align}
Thus, by taking the full expectation on both sides, we obtain the desired result (\ref{th3}).
\end{pf}

\begin{rmk}
\label{rmk6}
Let $u$ and $v$ are two unit vectors, i.e., $\|u\|_2=1$ and $\|v\|_2=1$. We use inner quantity $\langle u, v\rangle^2$ to represent the angle between $u$ and $v$, and the bigger the angle is, the bigger the fluctuation becomes from $u$ to $v$. Theorem \ref{theorem3} shows the fluctuation of two adjacent iterations. Specifically, when $\|A^T\frac{Ax_{k }- Ax_{\star}}{\|Ax_{k }- Ax_{\star}\|_2}\|_2^2$ is large, the angle between $\frac{Ax_{k}- Ax_{\star}}{\|Ax_{k}- Ax_{\star}\|_2}$ and $\frac{Ax_{k+1}- Ax_{\star}}{\|Ax_{k+1}- Ax_{\star}\|_2}$ is large, which implies that $\frac{Ax_{k}- Ax_{\star}}{\|Ax_{k}- Ax_{\star}\|_2}$ has a large fluctuation; when $\|A^T\frac{Ax_{k }- Ax_{\star}}{\|Ax_{k }- Ax_{\star}\|_2}\|_2^2$ is small, the angle between $\frac{Ax_{k}- Ax_{\star}}{\|Ax_{k}- Ax_{\star}\|_2}$ and $\frac{Ax_{k+1}- Ax_{\star}}{\|Ax_{k+1}- Ax_{\star}\|_2}$ is small, which implies that $\frac{Ax_{k}- Ax_{\star}}{\|Ax_{k}- Ax_{\star}\|_2}$ has very little fluctuation.

Since $\|A^T\frac{Ax_{k}- Ax_{\star}}{\|Ax_{k}- Ax_{\star}\|_2}\|_2^2\geq\sigma_r^2(A)$, Theorem \ref{theorem3} implies that if $Ax_{k}-Ax_{\star}$ is mainly composed of left singular vectors corresponding to the small singular values of $A$, its direction hardly changes, which means that the RGS method finally converges along left singular vector corresponding to the small singular value of $A$.

\end{rmk}

\section{Convergence directions of the REGS method}\label{sec4}
Recalling Remark \ref{rmk1}, when the matrix $A$ is not full column rank, the sequence $\{x_{k}\}_{k=1}^{\infty}$ generated by the RGS method does not converge to the minimum Euclidean norm least squares solution $x_{\star}$, even though $Ax_{k}$ does converge to $Ax_{\star}$. In \cite{Ma2015}, Ma et al. proposed an extended variant of the RGS method, i.e., the REGS method, to allow for convergence to $x_{\star}$ regardless of whether $A$ has full column rank or not.

Now, we list the REGS method presented in \cite{ Dukui2019} in Algorithm \ref{alg2}, which is a equivalent variant of the original REGS method \cite{Ma2015}, and restate its convergence bound presented in \cite{Dukui2019} in Theorem \ref{theorem5}. From the algorithm we find that, in each iteration, $x_k$ is the $k$th approximation of the RGS method and $z_k$ is a one-step RK update for the linear system $Az=Ax_{k}$ from $z_{k-1}$.

\begin{alg}
\label{alg2}
 The REGS method
\begin{enumerate}[]
\item \mbox{INPUT:} ~$A$, $b$, $\ell$, $x_{0}\in R^{n }$, $z_{0}\in \mathcal{R}(A^T)$
\item For $k=1, 2, \ldots, \ell-1$ do
\item ~~~~Select $j\in [n]$ with probability $\frac{\|A_{(j)}\|^2_2}{\|A\|^2_F}$
\item ~~~~Set $x_{k}=x_{k-1}-\frac{A_{(j)}^T ( Ax_{k-1}-b)}{ \| A_{(j)} \|_{2}^{2}}e_{(j)}$
\item ~~~~Select $i\in [m]$ with probability $\frac{\|A^{(i)}\|^2_2}{\|A\|^2_F}$
\item ~~~~Set $z_{k}=z_{k-1}-\frac{A^{(i)} ( z_{k-1}-x_{k})}{ \| A^{(i)} \|_{2}^{2}}(A^{(i)})^T$
\item End for
\end{enumerate}
\end{alg}

\begin{thm}
\label{theorem5}\cite{Dukui2019}
 Let $A\in R^{m\times n}$, $b\in R^{m}$, $x_{\star}=A^{\dag}b$ be the minimum Euclidean norm least squares solution, and $z_k$ be the $k$th approximation of the REGS method generated by Algorithm \ref{alg2} with initial $x_{0}\in R^{n }$ and $z_o\in \mathcal{R}(A^T)$. Then
\begin{align}
\mathbb{E}[\| z_{k}- x_{\star}\|_2^2]\leq(1-\frac{\sigma_r^2(A)}{\|A\|^2_F})^k \|z_{0}- x_{\star}\|_2^2+\frac{k}{\|A\|_F^2}(1-\frac{\sigma_r^2(A)}{\|A\|^2_F})^k \|Ax_0-Ax_{\star}\|_2^2.\label{th5}
\end{align}
\end{thm}

For the REGS method, we first discuss the convergence behavior of $z_{k}- x_{\star}$ in Theorem \ref{theorem6} and Theorem \ref{theorem7}, and then consider its convergence behavior of $Az_{k}- Ax_{\star}$ in Theorem \ref{theorem8}.

\begin{thm}
\label{theorem6}
 Let $A\in R^{m\times n}$, $b\in R^{m}$, $x_{\star}=A^{\dag}b$ be the minimum Euclidean norm least squares solution, and $z_k$ be the $k$th approximation of the REGS method generated by Algorithm \ref{alg2} with initial $x_{0}\in R^{n }$ and $z_o\in \mathcal{R}(A^T)$. Then
\begin{align}
\mathbb{E}[\langle z_{k}- x_{\star}, v_{\ell} \rangle]= (1 -\frac{\sigma_{\ell}^2(A)}{\|A\|^2_F} )^k    \langle z_{0}-x_{\star}, v_{\ell}\rangle  +\frac{k}{\|A\|^2_F} (1-\frac{\sigma_\ell^2(A)}{\|A\|^2_F})^k \langle A^T(Ax_{0}- Ax_{\star}),      v_{\ell} \rangle .\label{th6}
\end{align}
\end{thm}
\begin{pf}
From Algorithm \ref{alg2}, 
we have
\begin{align}
&\mathbb{E}[\langle z_{k}- x_{\star}, v_{\ell} \rangle] \notag
 \\
&= \mathbb{E}[\langle z_{k-1}-\frac{A^{(i)} ( z_{k-1}-x_{k})}{ \| A^{(i)} \|_{2}^{2}}(A^{(i)})^T- x_{\star}, v_{\ell} \rangle] \notag
 \\
&= \mathbb{E}[\langle (I-\frac{(A^{(i)})^TA^{(i)}}{\| A^{(i)} \|_{2}^{2}})( z_{k-1}-x_{\star}), v_{\ell}  \rangle] +\mathbb{E}[\langle \frac{(A^{(i)})^TA^{(i)}}{\| A^{(i)} \|_{2}^{2}}( x_{k}-x_{\star}),  v_{\ell}  \rangle],  \label{th65}
\end{align}
so we next consider $\mathbb{E}[\langle (I-\frac{(A^{(i)})^TA^{(i)}}{\| A^{(i)} \|_{2}^{2}})( z_{k-1}-x_{\star}), v_{\ell}  \rangle]$ and $\mathbb{E}[\langle \frac{(A^{(i)})^TA^{(i)}}{\| A^{(i)} \|_{2}^{2}}( x_{k}-x_{\star}),  v_{\ell}  \rangle]$ separately.

We first consider$\mathbb{E}[\langle (I-\frac{(A^{(i)})^TA^{(i)}}{\| A^{(i)} \|_{2}^{2}})( z_{k-1}-x_{\star}), v_{\ell}  \rangle]$. Let $\mathbb{E}_{k-1}[\cdot]$ be the conditional expectation conditioned on the first $k-1$ iterations of the REGS method. That is,
\begin{align}
 \mathbb{E}_{k-1}[\cdot]= \mathbb{E}[\cdot|j_1, i_1, j_2, i_2, \ldots, j_{k-1}, i_{k-1}],  \notag
\end{align}
where $j_{t^*}$ is the ${t^*}$th column chosen and $i_{t^*}$ is the ${t^*}$th row chosen. We denote the conditional expectation conditioned on
the first $k-1$ iterations and the $k$th column chosen as
\begin{align}
 \mathbb{E}_{k-1}^{i}[\cdot]= \mathbb{E}[\cdot|j_1, i_1, j_2, i_2, \ldots, j_{k-1}, i_{k-1}, j_k]. \notag
\end{align}
Similarly, we denote the conditional expectation conditioned on the first $k-1$ iterations and the $k$th row chosen as
\begin{align}
 \mathbb{E}_{k-1}^{j}[\cdot]= \mathbb{E}[\cdot|j_1, i_1, j_2, i_2, \ldots, j_{k-1}, i_{k-1}, i_k]. \notag
\end{align}
Then, by the law of total expectation, we have
\begin{align}
 \mathbb{E}_{k-1}[\cdot]=  \mathbb{E}_{k-1}^{j}[ \mathbb{E}_{k-1}^{i}[\cdot] ]. \notag
\end{align}
Thus, we obtain
\begin{align}
&\mathbb{E}_{k-1}[\langle (I-\frac{(A^{(i)})^TA^{(i)}}{\| A^{(i)} \|_{2}^{2}})( z_{k-1}-x_{\star}), v_{\ell}  \rangle] \notag
 \\
&=\mathbb{E}_{k-1}[\langle z_{k-1}-x_{\star}, v_{\ell}  \rangle -\langle \frac{A^{(i)} ( z_{k-1}-x_{\star})}{\| A^{(i)} \|_{2}^{2}}(A^{(i)})^T, v_{\ell}  \rangle] \notag
 \\
&=\langle z_{k-1}-x_{\star}, v_{\ell}  \rangle -\frac{1}{\|A\|^2_F} \sum\limits_{i=1}^{m}\langle A^{(i)} ( z_{k-1}-x_{\star})(A^{(i)})^T, v_{\ell}  \rangle  \notag
 \\
&=\langle z_{k-1}-x_{\star}, v_{\ell}  \rangle -\frac{1}{\|A\|^2_F} \sum\limits_{i=1}^{m}\langle(A^{(i)})^T,  z_{k-1}-x_{\star} \rangle\langle (A^{(i)})^T, v_{\ell}  \rangle  \notag
 \\
&=\langle z_{k-1}-x_{\star}, v_{\ell}  \rangle -\frac{1}{\|A\|^2_F}  \langle A ( z_{k-1}-x_{\star}),A v_{\ell}  \rangle.  \notag
\end{align}
Further, by making use of $z_{k-1}-x_{\star}=\sum\limits_{i=1}^{n}\langle z_{k-1}-x_{\star}, v_i\rangle v_i$ and $Av_i=\sigma_i(A) u_i $, we get
\begin{align}
&\mathbb{E}_{k-1}[\langle (I-\frac{(A^{(i)})^TA^{(i)}}{\| A^{(i)} \|_{2}^{2}})( z_{k-1}-x_{\star}), v_{\ell}  \rangle] \notag
\\
&=\langle z_{k-1}-x_{\star}, v_{\ell}  \rangle -\frac{1}{\|A\|^2_F}  \langle A  \sum\limits_{i=1}^{n}\langle z_{k-1}-x_{\star}, v_i\rangle v_i,\sigma_{\ell}(A) u_{\ell} \rangle   \notag
\\
&=\langle z_{k-1}-x_{\star}, v_{\ell}  \rangle -\frac{1}{\|A\|^2_F}  \langle   \sum\limits_{i=1}^{n}\langle z_{k-1}-x_{\star}, v_i\rangle \sigma_i(A) u_i,\sigma_{\ell}(A) u_{\ell} \rangle,  \notag
\end{align}
which together with the orthogonality of the left singular vectors $u_i$ yields
\begin{align}
&\mathbb{E}_{k-1}[\langle (I-\frac{(A^{(i)})^TA^{(i)}}{\| A^{(i)} \|_{2}^{2}})( z_{k-1}-x_{\star}), v_{\ell}  \rangle] \notag
\\
&=\langle z_{k-1}-x_{\star}, v_{\ell}  \rangle -\frac{\sigma_{\ell}^2(A)}{\|A\|^2_F}    \langle z_{k-1}-x_{\star}, v_{\ell}\rangle   \notag
\\
&=(1 -\frac{\sigma_{\ell}^2(A)}{\|A\|^2_F} )   \langle z_{k-1}-x_{\star}, v_{\ell}\rangle .  \notag
\end{align}
As a result, by taking the full expectation on both sides, we have
\begin{align}
 \mathbb{E}[\langle (I-\frac{(A^{(i)})^TA^{(i)}}{\| A^{(i)} \|_{2}^{2}})( z_{k-1}-x_{\star}), v_{\ell}  \rangle]=(1 -\frac{\sigma_{\ell}^2(A)}{\|A\|^2_F} ) \mathbb{E}[  \langle z_{k-1}-x_{\star}, v_{\ell}\rangle ].  \label{th62}
\end{align}

We now consider $\mathbb{E}[\langle \frac{(A^{(i)})^TA^{(i)}}{\| A^{(i)} \|_{2}^{2}}( x_{k}-x_{\star}),  v_{\ell}  \rangle]$. It follows from
\begin{align}
&\mathbb{E}_{k-1}[\langle \frac{(A^{(i)})^TA^{(i)}}{\| A^{(i)} \|_{2}^{2}}( x_{k}-x_{\star}),  v_{\ell}  \rangle] \notag
\\
&= \mathbb{E}_{k-1}^{j}[ \mathbb{E}_{k-1}^{i} [ \langle \frac{(A^{(i)})^TA^{(i)}}{\| A^{(i)} \|_{2}^{2}}( x_{k}-x_{\star}),  v_{\ell}  \rangle]]\notag
\\
&= \mathbb{E}_{k-1}^{j}[\frac{1}{\|A\|^2_F} \sum\limits_{i=1}^{m} \langle (A^{(i)})^TA^{(i)}( x_{k}-x_{\star}),  v_{\ell}  \rangle] \notag
\\
&= \mathbb{E}_{k-1}^{j}[\frac{1}{\|A\|^2_F} \sum\limits_{i=1}^{m} \langle  A^{(i)}( x_{k}-x_{\star}),  A^{(i)} v_{\ell}  \rangle] \notag
\\
&= \mathbb{E}_{k-1} [\frac{1}{\|A\|^2_F}   \langle  A ( x_{k}-x_{\star}),  A   v_{\ell}  \rangle], \notag
\end{align}
that
\begin{align}
 \mathbb{E} [\langle \frac{(A^{(i)})^TA^{(i)}}{\| A^{(i)} \|_{2}^{2}}( x_{k}-x_{\star}),  v_{\ell}  \rangle]=\mathbb{E}  [\frac{1}{\|A\|^2_F}   \langle  A ( x_{k}-x_{\star}),  A   v_{\ell}  \rangle].\label{th63}
\end{align}
Since
\begin{align}
 \mathbb{E}  [\frac{1}{\|A\|^2_F}   \langle  A ( x_{k}-x_{\star}),  A   v_{\ell}  \rangle]=\frac{\sigma_{\ell} (A)}{\|A\|^2_F} \mathbb{E}  [  \langle  A ( x_{k}-x_{\star}), u_{\ell}  \rangle], \notag
\end{align}
by exploiting (\ref{th1}) in Theorem \ref{theorem1}, we get
\begin{align}
&\mathbb{E}  [\frac{1}{\|A\|^2_F}   \langle  A ( x_{k}-x_{\star}),  A   v_{\ell}  \rangle] \notag
\\
&=\frac{\sigma_{\ell} (A)}{\|A\|^2_F} (1-\frac{\sigma_\ell^2(A)}{\|A\|^2_F})^k \langle Ax_{0}- Ax_{\star}, u_{\ell} \rangle \notag
\\
&=\frac{1}{\|A\|^2_F} (1-\frac{\sigma_\ell^2(A)}{\|A\|^2_F})^k \langle Ax_{0}- Ax_{\star},  A   v_{\ell} \rangle .\notag
\end{align}
Thus, substituting the above equality into (\ref{th63}), we have
\begin{align}
 \mathbb{E} [\langle \frac{(A^{(i)})^TA^{(i)}}{\| A^{(i)} \|_{2}^{2}}( x_{k}-x_{\star}),  v_{\ell}  \rangle]
 &=\frac{1}{\|A\|^2_F} (1-\frac{\sigma_\ell^2(A)}{\|A\|^2_F})^k \langle Ax_{0}- Ax_{\star},  A   v_{\ell} \rangle  \notag
 \\
 &=\frac{1}{\|A\|^2_F} (1-\frac{\sigma_\ell^2(A)}{\|A\|^2_F})^k \langle A^T(Ax_{0}- Ax_{\star}),      v_{\ell} \rangle . \label{th64}
\end{align}

Combining (\ref{th65}), (\ref{th62}) and (\ref{th64}) yields
\begin{align}
&\mathbb{E}[\langle z_{k}- x_{\star}, v_{\ell} \rangle] \notag
 \\
&= \mathbb{E}[\langle (I-\frac{(A^{(i)})^TA^{(i)}}{\| A^{(i)} \|_{2}^{2}})( z_{k-1}-x_{\star}), v_{\ell}  \rangle] +\mathbb{E}[\langle \frac{(A^{(i)})^TA^{(i)}}{\| A^{(i)} \|_{2}^{2}}( x_{k}-x_{\star}),  v_{\ell}  \rangle] \notag
 \\
&=(1 -\frac{\sigma_{\ell}^2(A)}{\|A\|^2_F} ) \mathbb{E}[  \langle z_{k-1}-x_{\star}, v_{\ell}\rangle ]+\frac{1}{\|A\|^2_F} (1-\frac{\sigma_\ell^2(A)}{\|A\|^2_F})^k \langle A^T(Ax_{0}- Ax_{\star}),      v_{\ell} \rangle \notag
 \\
&=(1 -\frac{\sigma_{\ell}^2(A)}{\|A\|^2_F} )^2 \mathbb{E}[  \langle z_{k-2}-x_{\star}, v_{\ell}\rangle ]+\frac{2}{\|A\|^2_F} (1-\frac{\sigma_\ell^2(A)}{\|A\|^2_F})^k \langle A^T(Ax_{0}- Ax_{\star}),      v_{\ell} \rangle \notag
 \\
&=\ldots=(1 -\frac{\sigma_{\ell}^2(A)}{\|A\|^2_F} )^k   \langle z_{0}-x_{\star}, v_{\ell}\rangle  +\frac{k}{\|A\|^2_F} (1-\frac{\sigma_\ell^2(A)}{\|A\|^2_F})^k \langle A^T(Ax_{0}- Ax_{\star}),      v_{\ell} \rangle, \notag
\end{align}
which is the desired result (\ref{th6}).

\end{pf}
\begin{rmk}
\label{rmk61}
 Theorem \ref{theorem6} shows that the decay rates of $\|z_k-x_{\star}\|_2$ are different in different right singular vectors spaces and the smallest singular value will lead to the slowest rate of convergence, which is the one in (\ref{th5}). So, the convergence bound presented by Du \cite{ Dukui2019} is optimal.

\end{rmk}

\begin{thm}
\label{theorem7}
 Let $A\in R^{m\times n}$, $b\in R^{m}$, $x_{\star}=A^{\dag}b$ be the minimum Euclidean norm least squares solution, and $z_k$ be the $k$th approximation of the REGS method generated by Algorithm \ref{alg2} with initial $x_{0}\in R^{n }$ and $z_o\in \mathcal{R}(A^T)$. Then
\begin{align}
\mathbb{E}[ \|z_{k}- x_{\star}\|^2_2]\leq \mathbb{E}[(1-\frac{1}{\|A\|^2_F}\|A\frac{z_{k-1}-x_{\star}}{\|z_{k-1}-x_{\star}\|_2}\|_2^2)\|z_{k-1}-x_{\star}\|_2^2]+ \frac{1}{\|A\|^2_F}(1-\frac{\sigma_r^2(A)}{\|A\|^2_F})^k \| Ax_{0}-Ax_{\star} \|^2_2. \label{th7}
\end{align}
\end{thm}

\begin{pf}
Following an analogous argument to Theorem 4 of \cite{Dukui2019}, we get
\begin{align}
\mathbb{E}   [\|z_{k}- x_{\star}\|_2^2 ]&= \mathbb{E}[ \|(I-\frac{(A^{(i)})^TA^{(i)}}{\| A^{(i)} \|_{2}^{2}})( z_{k-1}-x_{\star})\|_2^2]+\mathbb{E} [\|  \frac{(A^{(i)})^TA^{(i)}}{\| A^{(i)} \|_{2}^{2}}( x_{k}-x_{\star})\|_2^2], \notag
\end{align}
\begin{align}
\mathbb{E}[ \|(I-\frac{(A^{(i)})^TA^{(i)}}{\| A^{(i)} \|_{2}^{2}})( z_{k-1}-x_{\star})\|_2^2]
&=\mathbb{E} [(z_{k-1}-x_{\star})^T(I-\frac{A^TA}{\|A\|_F^2})(z_{k-1}-x_{\star})]\notag
\\
&=\mathbb{E}[(\|z_{k-1}-x_{\star}\|_2^2- \frac{1}{\|A\|^2_F}\|A(z_{k-1}-x_{\star})\|_2^2)]\notag
\\
&=\mathbb{E}[(1-\frac{1}{\|A\|^2_F}\|A\frac{z_{k-1}-x_{\star}}{\|z_{k-1}-x_{\star}\|_2}\|_2^2)\|z_{k-1}-x_{\star}\|_2^2],\notag
\end{align}
and
\begin{align}
\mathbb{E} [\|  \frac{(A^{(i)})^TA^{(i)}}{\| A^{(i)} \|_{2}^{2}}( x_{k}-x_{\star})\|_2^2]\leq\frac{1}{\|A\|^2_F}(1-\frac{\sigma_r^2(A)}{\|A\|^2_F})^k \| Ax_{0}-Ax_{\star} \|^2_2. \notag
\end{align}
Combining the above three equations, we have
\begin{align}
\mathbb{E}  [ \|z_{k}- x_{\star}\|_2^2]
\leq\mathbb{E}[(1-\frac{1}{\|A\|^2_F}\|A\frac{z_{k-1}-x_{\star}}{\|z_{k-1}-x_{\star}\|_2}\|_2^2)\|z_{k-1}-x_{\star}\|_2^2]+ \frac{1}{\|A\|^2_F}(1-\frac{\sigma_r^2(A)}{\|A\|^2_F})^k \| Ax_{0}-Ax_{\star} \|^2_2, \notag
\end{align}
which implies the desired result (\ref{th7}).
\end{pf}

\begin{rmk}
\label{rmk71}
Since $\|A\frac{z_{k-1}-x_{\star}}{\|z_{k-1}-x_{\star}\|_2}\|_2^2\geq\sigma_r^2(A)$, Theorem \ref{theorem7} implies that $z_{k}$ of the REGS method actually converges faster if $z_{k-1}-x_{\star}$ is not close to right singular vectors corresponding to the small singular values of $A$ .

\end{rmk}

\begin{thm}
\label{theorem8}
 Let $A\in R^{m\times n}$, $b\in R^{m}$, $x_{\star}=A^{\dag}b$ be the minimum Euclidean norm least squares solution, and $z_k$ be the $k$th approximation of the REGS method generated by Algorithm \ref{alg2} with initial $x_{0}\in R^{n }$ and $z_o\in \mathcal{R}(A^T)$. Then
\begin{align}
\mathbb{E}[\langle A z_{k}- Ax_{\star}, u_{\ell} \rangle]=(1 -\frac{\sigma_{\ell}^2(A)}{\|A\|^2_F} )^k    \langle Az_{0}-Ax_{\star}, v_{\ell}\rangle  +\frac{k}{\|A\|^2_F} (1-\frac{\sigma_\ell^2(A)}{\|A\|^2_F})^k \langle AA^T(Ax_{0}- Ax_{\star}),      u_{\ell} \rangle  .\label{th8}
\end{align}
\end{thm}

\begin{pf}
Similar to the proof of (\ref{th65}) in Theorem \ref{theorem6}, we obtain
\begin{align}
\mathbb{E}[\langle Az_{k}- Ax_{\star}, u_{\ell} \rangle]
 = \mathbb{E}[\langle A(I-\frac{(A^{(i)})^TA^{(i)}}{\| A^{(i)} \|_{2}^{2}})( z_{k-1}-x_{\star}), u_{\ell}  \rangle] +\mathbb{E}[\langle A \frac{(A^{(i)})^TA^{(i)}}{\| A^{(i)} \|_{2}^{2}}( x_{k}-x_{\star}),  u_{\ell}  \rangle].  \label{th80}
\end{align}
Then, we consider $\mathbb{E}[\langle A(I-\frac{(A^{(i)})^TA^{(i)}}{\| A^{(i)} \|_{2}^{2}})( z_{k-1}-x_{\star}), u_{\ell} \rangle]$ and $\mathbb{E}[\langle  A \frac{(A^{(i)})^TA^{(i)}}{\| A^{(i)} \|_{2}^{2}}( x_{k}-x_{\star}),  u_{\ell}   \rangle]$ separately.

We first consider $\mathbb{E}[\langle A(I-\frac{(A^{(i)})^TA^{(i)}}{\| A^{(i)} \|_{2}^{2}})( z_{k-1}-x_{\star}), u_{\ell} \rangle]$. It follows from
$$ \langle A(I-\frac{(A^{(i)})^TA^{(i)}}{\| A^{(i)} \|_{2}^{2}})( z_{k-1}-x_{\star}), u_{\ell} \rangle = \langle (I-\frac{(A^{(i)})^TA^{(i)}}{\| A^{(i)} \|_{2}^{2}})( z_{k-1}-x_{\star}), A^Tu_{\ell} \rangle $$ and $A^Tu_{\ell}=\sigma_{\ell}(A)v_{\ell}$, that
\begin{align}
  \mathbb{E}[\langle A(I-\frac{(A^{(i)})^TA^{(i)}}{\| A^{(i)} \|_{2}^{2}})( z_{k-1}-x_{\star}), u_{\ell} \rangle] = \sigma_{\ell}(A)\mathbb{E}[\langle  (I-\frac{(A^{(i)})^TA^{(i)}}{\| A^{(i)} \|_{2}^{2}})( z_{k-1}-x_{\star}),  v_{\ell} \rangle],  \notag
\end{align}
which together with (\ref{th62}), yields
\begin{align}
  \mathbb{E}[\langle A(I-\frac{(A^{(i)})^TA^{(i)}}{\| A^{(i)} \|_{2}^{2}})( z_{k-1}-x_{\star}), u_{\ell} \rangle]
  &= \sigma_{\ell}(A)(1 -\frac{\sigma_{\ell}^2(A)}{\|A\|^2_F} ) \mathbb{E}[  \langle z_{k-1}-x_{\star}, v_{\ell}\rangle ]  \notag
  \\
  &=(1 -\frac{\sigma_{\ell}^2(A)}{\|A\|^2_F} ) \mathbb{E}[  \langle z_{k-1}-x_{\star},\sigma_{\ell}(A) v_{\ell}\rangle ]\notag
    \\
  &=(1 -\frac{\sigma_{\ell}^2(A)}{\|A\|^2_F} ) \mathbb{E}[  \langle Az_{k-1}-Ax_{\star},u_{\ell}\rangle ].  \label{th81}
\end{align}

We now consider $\mathbb{E}[\langle  A \frac{(A^{(i)})^TA^{(i)}}{\| A^{(i)} \|_{2}^{2}}( x_{k}-x_{\star}),  u_{\ell}   \rangle]$. Exploiting (\ref{th64}), we have
\begin{align}
 \mathbb{E}[\langle  A \frac{(A^{(i)})^TA^{(i)}}{\| A^{(i)} \|_{2}^{2}}( x_{k}-x_{\star}),  u_{\ell}   \rangle]
&=\mathbb{E}[\langle    \frac{(A^{(i)})^TA^{(i)}}{\| A^{(i)} \|_{2}^{2}}( x_{k}-x_{\star}), A^T u_{\ell}   \rangle]\notag
\\
&=\sigma_\ell (A)\mathbb{E}[\langle    \frac{(A^{(i)})^TA^{(i)}}{\| A^{(i)} \|_{2}^{2}}( x_{k}-x_{\star}), v_{\ell}   \rangle]\notag
 \\
 &=\frac{\sigma_\ell (A)}{\|A\|^2_F} (1-\frac{\sigma_\ell^2(A)}{\|A\|^2_F})^k \langle A^T(Ax_{0}- Ax_{\star}),      v_{\ell} \rangle \notag
  \\
 &=\frac{1}{\|A\|^2_F} (1-\frac{\sigma_\ell^2(A)}{\|A\|^2_F})^k \langle AA^T(Ax_{0}- Ax_{\star}),      u_{\ell} \rangle. \label{th82}
\end{align}

Thus, combining (\ref{th80}), (\ref{th81}) and (\ref{th82}) yields
\begin{align}
&\mathbb{E}[\langle Az_{k}- Ax_{\star}, u_{\ell} \rangle]  \notag
 \\
&= \mathbb{E}[\langle A(I-\frac{(A^{(i)})^TA^{(i)}}{\| A^{(i)} \|_{2}^{2}})( z_{k-1}-x_{\star}), u_{\ell}  \rangle] +\mathbb{E}[\langle A \frac{(A^{(i)})^TA^{(i)}}{\| A^{(i)} \|_{2}^{2}}( x_{k}-x_{\star}),  u_{\ell}  \rangle] \notag
 \\
&=(1 -\frac{\sigma_{\ell}^2(A)}{\|A\|^2_F} ) \mathbb{E}[  \langle Az_{k-1}-Ax_{\star},u_{\ell}\rangle ]+\frac{1}{\|A\|^2_F} (1-\frac{\sigma_\ell^2(A)}{\|A\|^2_F})^k \langle AA^T(Ax_{0}- Ax_{\star}),      u_{\ell} \rangle \notag
 \\
&=(1 -\frac{\sigma_{\ell}^2(A)}{\|A\|^2_F} )^2 \mathbb{E}[  \langle Az_{k-2}-Ax_{\star},u_{\ell}\rangle ]+\frac{2}{\|A\|^2_F} (1-\frac{\sigma_\ell^2(A)}{\|A\|^2_F})^k \langle AA^T(Ax_{0}- Ax_{\star}),      u_{\ell} \rangle \notag
 \\
&=\ldots=(1 -\frac{\sigma_{\ell}^2(A)}{\|A\|^2_F} )^k   \langle Az_{0}-Ax_{\star}, v_{\ell}\rangle  +\frac{k}{\|A\|^2_F} (1-\frac{\sigma_\ell^2(A)}{\|A\|^2_F})^k \langle AA^T(Ax_{0}- Ax_{\star}),      u_{\ell} \rangle, \notag
\end{align}
which is the desired result (\ref{th8}).
\end{pf}

\begin{rmk}
\label{rmk81}

Theorem \ref{theorem8} shows the decay rates of $\|Az_{k}- Ax_{\star}\|$ of the REGS method and suggests that small singular values lead to poor convergence rates and vice versa. We note similar issues arise for the RK, REK, and RGS methods discussed in \cite{steinerberger2021randomized}, \cite{zhang2021preconvergence}, and Theorem \ref{theorem1}, respectively.

\end{rmk}

\section{Numerical experiments }\label{sec5}
Now we present two simple examples to illustrate that the convergence directions of the RGS and REGS methods. To this end, let $G_0\in R^{500\times 500}$ be a Gaussian matrix with i.i.d. $N(0, 1)$ entries and $D\in R^{500\times 500}$ be a diagonal matrix whose diagonal elements are all 100. Further, we set $G_1=G_0+D$ and replace 
its last row $G_1^{(500)}$ by a tiny perturbation of $G_1^{(499)}$, i.e., adding 0.01 to each entry of $G_1^{(499)}$. Then, we normalize all rows of $G_1$, i.e., set $\|G_1^{(i)}\|_2=1$, $i=1, 2, \ldots, 500$. After that, we set
$A_1=\begin{bmatrix}
G_1\\
G_2
\end{bmatrix}
\in R^{600\times 500}$ and
$A_2=\begin{bmatrix}
G_1,
G_3
\end{bmatrix}
\in R^{500\times 600}$, where $G_2\in R^{100\times 500}$ and $G_3\in R^{500\times 100}$ are zero matrices. So, the first 499 singular values of the matrices $A_1$ and $A_2$ are between $\sim 0.6$ and $\sim 1.5$, and the smallest nonzero singular value is $\sim 10^{-4}$.

\begin{figure}[ht]
 \begin{center}
\includegraphics [height=5.5cm,width=8.5cm  ]{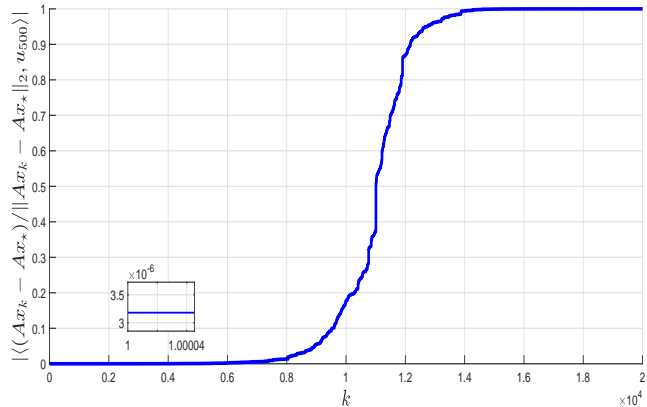}
 \end{center}
\caption{A sample evolution of $ |\langle (Ax_k-Ax_{\star})/\|Ax_k-Ax_{\star}\|_2, u_{500} \rangle|$ of the RGS method. }\label{fig1}
\end{figure}

\begin{figure}[ht]
 \begin{center}
\includegraphics [height=5.5cm,width=8.5cm  ]{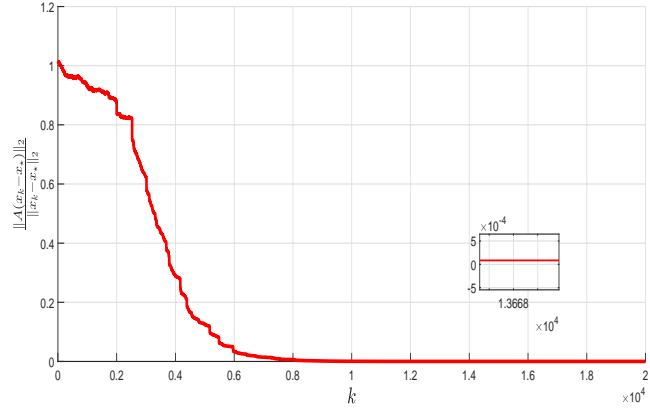}
 \end{center}
\caption{A sample evolution of $\frac{\|A(x_k-x_{\star})\|_2 }{\|x_k-x_{\star}\|_2}$ of the RGS method. }\label{fig2}
\end{figure}

We first consider convergence directions of $Ax_k-Ax_{\star}$ and $x_k-x_{\star}$ of the RGS method. We generate a vector $x \in R^{500}$ using the MATLAB function \texttt{randn}, set the full column rank coefficient matrix $A=A_1$ and set the right-hand side $b=A x +z$, where $z$ is a nonzero vector belonging to the null space of $A^{T}$, which is generated by the MATLAB function \texttt{null}. With $x_0=0$, we plot $|\langle (Ax_k-Ax_{\star})/\|Ax_k-Ax_{\star}\|_2, u_{500} \rangle|$ and $\frac{\|A(x_k-x_{\star})\|_2 }{\|x_k-x_{\star}\|_2}$ in Figure \ref{fig1} and Figure \ref{fig2}, respectively.

From Figure \ref{fig1}, we find that $|\langle (Ax_k-Ax_{\star})/\|Ax_k-Ax_{\star}\|_2, u_{500} \rangle|$ initially is very small and almost is 0, which indicates that $Ax_k-Ax_{\star} $ is not close to the left singular vector $u_{500}$. Considering the analysis of Remark \ref{rmk5}, the phenomenon implies the `preconvergence' behavior of the RGS method, that is, the RGS method seems to converge quickly at the beginning. In addition, as $k\rightarrow\infty$, $|\langle (Ax_k-Ax_{\star})/\|Ax_k-Ax_{\star}\|_2, u_{500} \rangle|\rightarrow 1$. This phenomenon implies that $Ax_{k}-Ax_{\star}$ tends to the left singular vector corresponding to the smallest singular value of $A$.

From Figure \ref{fig2}, we observe that the values of $\frac{\|A(x_k-x_{\star})\|_2 }{\|x_k-x_{\star}\|_2}$ decreases with $k$ and finally approaches the small singular value. This phenomenon implies that the forward direction of $x_k-x_{\star}$ is mainly determined by the right singular vectors corresponding to the large singular values of $A$ at the beginning. With the increase of $k$, the direction is mainly determined by the right singular vectors corresponding to the small singular values. Finally, $x_k-x_{\star}$ tends to the right singular vector space corresponding to the smallest singular value. Furthermore, this phenomenon also allows for an interesting application, i.e., finding nonzero vectors $x$ such that $\frac{\|Ax\|_2}{\|x\|_2}$ is small.

\begin{figure}[ht]
 \begin{center}
\includegraphics [height=5.5cm,width=8.5cm  ]{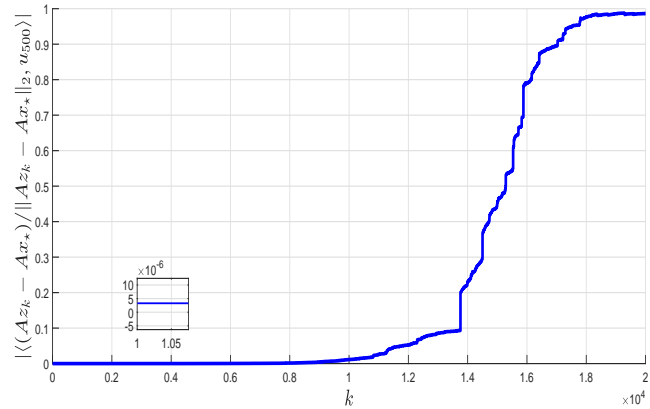}
 \end{center}
\caption{A sample evolution of $ |\langle (Az_k-Ax_{\star})/\|Az_k-Ax_{\star}\|_2, u_{500} \rangle|$ of the REGS method. }\label{fig3}
\end{figure}

\begin{figure}[ht]
 \begin{center}
\includegraphics [height=5.5cm,width=8.5cm  ]{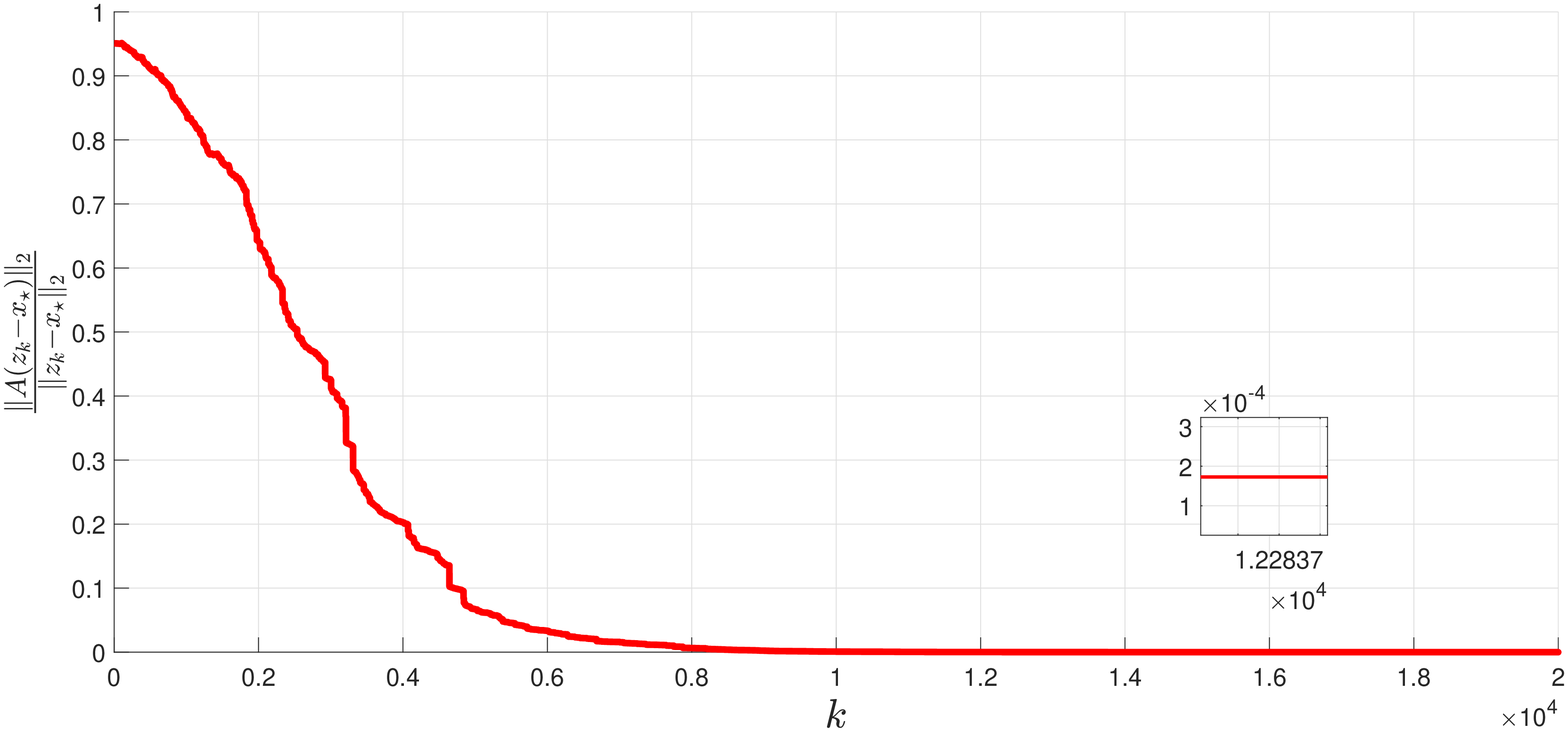}
 \end{center}
\caption{A sample evolution of $\frac{\|A(z_k-x_{\star})\|_2 }{\|z_k-x_{\star}\|_2}$ of the REGS method. }\label{fig4}
\end{figure}

We now consider convergence directions of $Az_k-Ax_{\star}$ and $z_k-x_{\star}$ of the REGS method. We generate a vector $x \in R^{600}$ using the MATLAB function \texttt{randn}, set the coefficient matrix $A=A_2$ which does not have full column rank, and set the right-hand side $b=Ax  $. With $x_0=0$ and $z_0=0$, we plot $|\langle (Az_k-Ax_{\star})/\|Az_k-Ax_{\star}\|_2, u_{500} \rangle|$ and $\frac{\|A(z_k-x_{\star})\|_2 }{\|z_k-x_{\star}\|_2}$ in Figure \ref{fig3} and Figure \ref{fig4}, respectively.

Figure \ref{fig3} and Figure \ref{fig4} show the similar results obtained in the RGS method. That is, the convergence directions of $A z_k-Ax_{\star} $ and $ z_k- x_{\star} $ of the REGS method initially are depending on the large singular values and then mainly depending on the small singular values, and finally depending on the smallest singular value of $A$.

\bibliography{mybibfile}

\end{document}